\documentclass[12pt]{article}
\usepackage[utf8]{inputenc}
\usepackage[russian]{babel}
\usepackage{amsthm}
\usepackage{amsmath}
\usepackage{amsfonts}
\usepackage{graphicx}
\usepackage{hyperref}
\usepackage{mathtools}
\usepackage{amssymb}
\usepackage{amsthm}
\usepackage[a4paper, left=10mm, right=16mm, top=20mm, bottom=20mm]{geometry}
\usepackage{amsfonts,amssymb,mathrsfs,amscd}
\def \R {{\mathbb {R}}}

\def \Z {{\mathbb {Z}}}

\def \P {\mathcal {P}}
\def \B {{\cal B}}
\def \M {{\mathcal  M}}

\title{\bf Asymmetry of Entropy Invariants \\
for Generic Mixing $\Z^n$-Actions}
\author{\Large \it  Mikhail V. Engelgardt, \ Valery V. Ryzhikov}
\date{}

\begin{document}
\maketitle
\large

\section{Introduction}
The Kirillov-Kushnirenko entropy (see \cite{Ki}--\cite{RT}) is a collective name for the continuum of $\P$-invariants for actions preserving the probability measure. Here $\P=\{P_j\}$ denotes a sequence of finite subsets $P_j$ of a countable group (we restrict our consideration to $\Z^n$-actions). If $P_j\subset \Z$ have the same finite cardinality, and the distances between elements of $P_j$ increase with $j$, then the corresponding $\P$-entropy is directly related to the property of multiple mixing, introduced by Rokhlin in \cite{Ro} to distinguish ordinary mixing from multiple mixing. Rokhlin sought corresponding examples among automorphisms of commutative  compact groups. However, it turned out that in this class of algebraic systems, ergodic $\Z$-actions possess mixing of all multiplicities.

Mixing $\Z$-actions that do not exhibit multiple mixing have not yet been found (for Rokhlin's problem and its relationship with other invariants of dynamical systems, see \cite{Ry}). Various entropy-type properties could have emerged before introduction of Kolmogorov's 
entropy  as a strengthening  of multiple mixing, if the mixing multiplicity is allowed to increase slowly. Invariants playing a similar role arose in the works of Kirillov \cite{Ki} and Kushnirenko \cite{Ku}, somewhat belatedly, as a generalization and reaction to the now fundamental concept of entropy per unit time. Kirillov defined them for group actions, and Kushnirenko provided the finiteness of the $2^n$-entropy of the automorphisms  $T$ entering the horocyclic flows. 
Thus, in the class of systems with classical zero entropy, he discovered nonisomorphic systems with the same Lebesgue spectrum.
These are the powers $T \times\dots\times T$, since their $2^n$-entropy values  are different.

Rokhlin emphasized the search, which was relevant at the dawn of ergodic theory,
for actions with the same spectrum but metrically (from the word "measure") non-isomorphic. Soon, Anzai \cite{Anz} introduced the concept of  skew product and, in particular, found a solution to this problem. He proposed  a skew product $T$, which is not isomorphic to $T^{-1}$.  However,
a commutative action and its inverse are spectrally isomorphic.
Later, Kolmogorov \cite{Ko}  introduced the  entropy, obtaining a
continuum of measure nonisomorphic transformations with the same Lebesgue
spectrum. Entropy is a symmetric numerical invariant, since its values  for the transformations and its inverse are the same. We note that one of the goals of this note was to find some asymmetric entropy invariant. For $\Z$-actions, unlike for the actions of $\Z^n$ groups for $n>1$, this problem has not yet been solved. There is hope to piece together the asymmetry bit by bit, as was done in \cite{aa} for other purposes. To do this, we need to find an entropy invariant such that, for a suitable action, the entropy characteristics along the sequence $\P$ gradually decrease to zero, i.e., the $\P$-invariant is zero, while for a symmetric sequence, $\P^\ast=-\P$ is positive.

Another important event related to Rokhlin's problem, the entropy, and asymmetry of action, occurred even later. Rokhlin's idea to use automorphisms of commutative groups to solve the multiple mixing problem is successful in the case of Z-actions for n>1. Corresponding examples were discovered by Ledrappier. His paper \cite{L} gives an elegant example of a symmetric action, but this example is obviously associated with an infinite class of actions, including asymmetric actions.
All these actions are called Ledrappier systems. We give the simplest example, which we will be used later.

Consider the group $X\subset \Z_2^{\Z^2}$ consisting of all sequences $x$ satisfying the identity
$$x+T^{(1,0)}x + T^{(0,1)}x =0. \eqno (1) $$
The action $\{T^z\}$ is defined by the equality $T^zx(v)=x(v-z)$. It preserves the Haar measure $\mu$, and
the group automorphisms of $T^z$ thus become automorphisms of the space $(X,\mu)$.
The action is asymmetric in the sense that identity (1) is violated for its inverse. We will exploit this and establish that for a specially chosen sequence $\P$ of finite subsets of the group $\Z^2$, our action has zero $\P$-entropy,
and its inverse has infinite completely positive $\P$-entropy. This implies not only non-isomorphism, but also a stronger property of disjointness of the action and its inverse.

Thus, multiple mixing and its absence, entropy, and the asymmetry of action in its strongest manifestation are intertwined in the class of Rokhlin-Ledrapier systems. To this, we add the typicality of such effects in the Alpern-Tikhonov sense. We will show that
the asymmetry of the $\P$-entropy for some sequences $\P$ is typical in the space of mixing $\Z^n$-actions for $n>1$. This implies the typicality of the $\Z^n$-action and its inverse for $n>0$. Here, we apply results of Bashtanov \cite{B}, Tikhonov \cite{T}, our theorem on entropy asymmetry, and the $\P$-analog of Pinsker theorem from the Ryzhikov-Thouvenot paper, \cite{RT}, Theorem 4.1.

In the case $n=1$, to prove the genericity of the strong asymmetry, we use Rudolph's example \cite{Ru} of a mixing automorphism disjoint with its powers. To our knowledge, this is the only example in the literature. It is possible that the Poisson suspensions from \cite{aa} are also suitable systems, but this has not been proven.

\section{Asymmetric $\P$-entropy}
Let $\P=\{P_j\}$ of finite sets $P_j\subset G$ be given, where $G$ is an infinite group, and let $\Phi=\{T^g\}$
be a measure-preserving $G$-action.
For a finite measurable partition $\xi=\{C_1,\dots,C_p\}$ of the space $X$, we define the quantities $$h_j(\Phi,\xi)=\frac 1 {|P_j|} H\left(\bigvee_{g\in P_j}T^g\xi\right),$$
where $H(\xi)=-\sum_{i=1}^p \mu( C_i)\ln \mu( C_i)$ is the entropy of the partition $\xi$.

\vspace{3mm}
\bf $\P$-entropy. \rm We define for an action $\Phi$ its entropy with respect to a partition $\xi$ as the quantity
$h_{P}(\Phi,\xi)={\limsup_j} \ h_k(\Phi,\xi).$
The $\P$-entropy of the action $\Phi$ is defined as
$h_{P}(\Phi)=\sup_\xi h_{P}(\Phi,\xi),$
the upper bound is taken over finite $\mu$-measurable partitions of $\xi$.

\vspace{3mm}
\bf  Ledrappier actions. \rm We have already described the $\Z^2$-action above. Consider the case $n=3$. The group $\Z_2^{\Z^3}$ contains a compact subgroup
$X$ formed by all sequences $x: \Z^3\to\Z_2$ satisfying the identity
$$x(z+(1,0,0))+x(z+(0,1,0)) +x(z+(0,0,1)) =x(z).$$
The action of the group $\Z^3$ is defined by the equality
$T^wx(z)=x(z-w).$ It is easy to verify that the induced identities
$$x(z+(2^n,0,0))+x(z+(0,2^n,0)) +x(z+(0,0,2^n)) =x(z).$$
This remarkable 2-similarity plays a key role in what follows.
The required $\Z^n$-actions for $n>3$ are defined similarly.
From now on, we will only consider Ledrappier's $\Z^2$-actions, since other cases are similar, but the formulas are longer.

\vspace{3mm}
\bf Chois of $\P$. \rm
We define $P_j\subset \Z^2$ by induction:
$$P_0=\{(0,0), (0,1), (1,0)\},$$
$$P_{j+1}= 2^{j}P_j\ \cup \ 2^{j}P_j+( 4^{j}, 0)\ \cup \ 2^{j}P_j+(0, 4^{j}).$$
Let $\xi=\{A,\, X\setminus A\}$, where $A=\{x\in X\, :\, x(0,0)=0\}.$

Denote $P_j^\ast=-P_j$.

\vspace{3mm}
\bf Theorem 2.1. \it $\P$-entropy of the action $\Phi$ is zero, and  $\P$-entropy of the action $\Phi^\ast$ is completely positive.\rm

\vspace{3mm}
The theorem follows directly from the following lemmas.

\vspace{3mm}
\bf Lemma 2.2. \it Let the partition $\nu$ be subordinate to a finite sum of shifts of the partition $\xi$. For all sufficiently large $j$, we have
$H(\vee_{z\in P^\ast_j} T_z \nu)=3^j H(\nu)$.\rm

\vspace{3mm}
In the space $L_2(X,\mu)$, the characters form an
orthonormal basis. Consider the characters $\chi_z$ of the group $X$ defined by the equalities
$$\chi_z(x)=1, \ x(z)=0; \ \ \ \chi_z(x)=-1, \ x(z)=1.$$
Let $\xi=\{A,\, X\setminus A\}$, where $A=\{x\in X\, :\, x(0,0)=0\}.$
The group of characters corresponds to a minimal algebra of sets with respect to which these characters are measurable. This algebra is generated by atoms of the same measure,
which is easily proved by induction. We only need to take into account that different characters are orthogonal and that all characters of the form
$\chi_{(q,0)}$ are independent: all possible products of them are nontrivial characters. The set algebra associated with the group generated by the characters $\chi_z$, $z\in P^\ast_j$,
contains $3^j$ atoms, which are the atoms of the partition $\xi_j=\bigvee_{z\in P_j} T^z\xi$. Shifting the algebra under the automorphism yields a similar shift of characters. But the set $P_j$ and its shifts are far apart for large $j$. The character group, its horizontal shift, and its vertical shift are independent, so the corresponding partitions
$$\xi_j, \ \ T^{(0,4^j)}\xi_j, \ \ T^{(0,4^j)}\xi_j.$$
All this leads to the assertion of the lemma (we leave the details to the interested reader).

\vspace{3mm}
\bf Lemma 2.3. \it $h_j(\Phi,\xi)= (2/3)^j \ln 2.$
\rm

\vspace{3mm}
Proof.
Note that the group $G_{P_j}$ generated by definition by all characters $\chi_z$ for $z\in P_j$ is isomorphic to the group
$G_{Q_j}$, where $Q_j= P_j\cap \{ (q,0): q\in \Z
\}$, since $$\chi_{(q,4^m)}= \chi_{(q,0)}\chi_{(q,4^m)}.$$
All characters in our group over the line $D=\{(w,0): w\in \Z\}$ are products of characters of the form $\chi_d, d\in D$.

The number of atoms of the partition $\vee_{z\in P^\ast_j} T_z \xi$ is $2^n$,
consequently, the action has zero $\P$-entropy.

\section{  Asymmetry is generic}
A number of properties of mixing $\Z^n$-actions of groups $\Z^n$ on the probability space $(X,\B,\mu)$ are typical. These include the multiple mixing property, established by Tikhonov (see \cite{T} and references), and, for example, the triviality of the action centralizer, shown by Bashtanov \cite{B}. A property being typical means that all actions with this property contain a dense set of type $G_\delta$. The topology in the space of mixing $\Z^n$-actions is determined by the metric $r_n$ defined below.

The group $Aut(\mu)$ of atomorphisms of the space $(X,\mu)$ is equipped with the Halmos metric $\rho$: for $S,T$:
$$ \rho(S,T)=\sum_i 2^{-i}\left(\mu(S^A_i\Delta TA_i)+\mu(S^{-1}A_i\Delta T^{-1}A_i)\right),$$
where $\{A_i\}$ is a fixed set dense in the algebra $\B$.
By summing the distances between the generators of $\Z^n$-actions,
we obtain the metric $\rho_n$ on the space of $\Z^n$-actions.

We introduce a metric $d$ on $Aut(\mu)$:
$$ d(S,T)=\sum_{i,j} 2^{-i-j}\left|\mu(SA_i\cap A_j)-\mu(TA_i\cap A_j)\right|.$$
On the set $Mix$ of all mixing $\Z^n$-actions, the metric $r_n$ is defined as follows:
$$ r_n(\{S^z\},\{T^z\})= \rho_n(\{S^z\},\{T^z\}) + \sup_{z\in \Z^n}d(S^z,T^z).$$
Recall that an action $\Phi=\{T^z\,:\, z\in\Z^n\}$ is called mixing if for every $A,B\in\B$, we have
$$ \mu(T^zA\cap B)\to \mu(A)\mu(B), \ z\to\infty.$$

For mixing $\Z^n$-actions, it is known (and we will re-prove it in a special case) that the infinite Kirillov-Kushnirenko entropy is typical for sequences of expanding sets $P_j$.

Consider the space $\M$ of measures on $X\times X$,
satisfying the following projection conditions: for every $A\in \B$, we have
$\nu(A\times X)=\nu(X\times A)=\mu(A).$
The set $\M$ is equipped with the metric $dist$:
$$ dist(\nu,\eta)=\sum_{i,j} 2^{-i-j}\left|\nu(A_i\times A_j)-\eta(A_i\times A_j)
\right|,$$ where $\{A_i\}$ is a fixed family dense in the algebra of all $\mu$-measurable sets. The space $(\M,dist)$ is compact.

\it The disjointness of the \rm automorphisms $S,T$ means that for all measures except $\mu^2=\mu\times \mu$, $dist(\nu, (S\times T)\nu)>0$.

\vspace{3mm}
\bf Theorem 3.1. \it A generic mixing $\Z^n$-action and its inverse are disjoint.\rm

\vspace{3mm}
Proof. The theorem for $n>1$ is a consequence of Theorem 2.1 and Theorem 4.1 from \cite{RT} and the following facts. It is easy to verify that actions with the properties asserted by Theorem 2.1 form a $G_\delta$-set. The conjugacy class of a mixing action is dense in the space of mixing actions, as shown by Bashtanov \cite{B}. One can also use results of Tikhonov \cite{T}, including the density of the conjugacy class of the Cartesian square of a mixing $\Z^n$-action (and $\R^n$-action).
The square inherits zero entropy of the action while simultaneously having a completely positive entropy of the inverse action. This ensures the density of the specified $G_\delta$-set. Note that this approach is effective in proving an analog of Theorem 3.1 for $\R^n$-actions (for which Bashtanov's theorem does not exist).

It remains to consider the case $n=1$. We will prove the following assertion and thus Theorem 3.1.

\vspace{3mm}
\bf Theorem 3.2. \it For a generic mixing automorphism $T$ of a probability
space, its degrees $T^p$ and $T^q$ for $p\neq q$ are disjoint.\rm

\vspace{3mm}
Proof. An example of a transformation with this property is given in Rudolf's paper \cite{Ru}.
The conjugacy class of every mixing automorphism is dense in $(Mix, r_1)$, as Bashtanov showed \cite{B}. Let us verify that all automorphisms with the specified property form a $G_\delta$-set. We consider only the case $p=1, q=-1$, since the proof of the other cases is completely similar. Set
$$K_m=\{\nu\in \M\ :\ dist(\nu,\mu^2)\ \geq \ 1/m\}.$$
The sets $K_m$ are closed in $(\M,dist)$, hence they are compact. For a fixed $T\in Aut$,
disjoint from $T^{-1}$, we have
$$ \min \{ dist(\nu\, ,\, (T\times T^{-1})\nu) \ :\ \nu\in K_m\}=d_m>0,$$
since a continuous positive function on the compact set $K_m$ is separated from 0.
We define the set
$$U_{m}= \bigcup_n \{ S\,:\, dist(\nu, \, (S\times S^{-1})\nu)> \ 1/n , \ \nu\in K_m\}.$$
It is open as the union of open sets.
It remains to note that $\bigcap_m U_m$, the family of all mixing automorphisms disjoint from their inverses, is a dense $G_\delta$-set. The theorem is proved.

\vspace{2mm}
\bf Theorem 3.3. \it Let $a$ be an endomorphism of the group $\Z^n$. If there exists a mixing action $\Phi\in Mix_n$ disjoint from $\tilde\Phi^a$, then the set of all mixing actions $\Phi$ disjoint from $\Phi^a$ is generic.\rm

\vspace{2mm}
This conditional statement is proved by the same method as Theorem 3.1. The density
of the conjugacy class of a mixing $\Z^n$-action is established in \cite{B}.
However, the required examples of actions are missing in the literature. Note that they can be constructed, but this requires considerable effort. It is possible that Ledrappier's actions provide the necessary examples.

In conclusion, let's note a small quasi-paradox. We obtained an asymmetry for typical actions from the absence
of multiple mixing. But typical mixing $\Z^n$-actions, as Tikhonov showed long ago, have mixing of all orders.

\section{ Remarks}
 \rm  We obtained an asymmetry for generic actions from the absence of
multiple mixing. But generic mixing $\Z^n$-actions, as Tikhonov showed \cite{T12}, have mixing of all multiplicities.

Let's  explain  the title of our paper. From Theorem 3.1.1 \cite{R21} (remark of V.R.: the word "amenable" was accidentally omitted in the formulation of this theorem)
it follows that for a generic $\Z^n$-action, its $P_j$ and $P_j^\ast$ entropies are infinite. The proof of genericity for mixing $\Z^n$-actions is exactly the same. Thus, both entropies are the same in our case. How does the asymmetry manifest itself? The answer is given by the following theorem, which we have conceptually proved, so we will limit ourselves to the formulation.

\vspace{3mm}
\bf Theorem 4.1. \it For our sequence $\P=\{P_j\}$, there is typically a set of mixing actions $\Phi$ for which there exists a proper sequence $j(k)\to\infty$ such that the $\{P_{j(k)}\}$-entropy of the action $\Phi$ is zero, and the $\{P^\ast_{j(k)}\}$-entropy is completely positive. In this case, for the inverse action $\Phi^\ast$
it is the opposite: the $\{P^\ast_{j(k)}\}$-entropy is zero, and the 
$\{P_{j(k)}\}$-entropy is completely positive. Moreover, the $P_j$ and $P_j^\ast$ entropies of these actions are infinite.\rm

\vspace{3mm}
Thanks to $\P$-entropy, we can enhance the effect of Kolmogorov entropy, which distinguishes a continuum of actions with a Lebesgue spectrum. One can specify a continuum of pairwise disjoint actions with a Lebesgue spectrum \cite{faa}. But these actions are symmetric, since they are factors of always symmetric Gaussian systems.

\large
\newpage
\begin{center}
{\bf \Large  Асимметрия энтропийных инвариантов \\
       для типичных перемешивающих $\Z^n$-действий }

\vspace{4mm}
{ \it \Large Валерий В. Рыжиков,\  Михаил В. Энгельгардт}

\end{center}

\vspace{4mm} 
{\bf 1. Введение}

\vspace{5mm}
Энтропия Кириллова-Кушниренко (см.\cite{Ki}--\cite{RT}) является  собирательным названием   континуума $\P$-инвариантов для  действий, сохраняющих вероятностную меру. Здесь  $\P=\{P_j\}$ обозначает последовательность конечных подмножеств $P_j$ счетной группы (мы ограничимся рассмотрением действий групп $\Z^n$).  Если множества $P_j\subset \Z$  имеют одинаковую конечную мощность, а расстояния между элементами множества $P_j$ увеличивается с ростом $j$,  то соответствующая $\P$-энтропия  непосредственно связана  со свойством кратного перемешивания, введенного Рохлиным в \cite{Ro}.   Чтобы отличить обычное перемешивание от кратного, Рохлин искал соответствующие примеры среди  автоморфизмов  компактных коммутативных  групп. Но оказалось, что в этом классе  алгебраических систем эргодические $\Z$-действия  обладают перемешиванием всех кратностей. 

Рохлин  акцентировал внимание на актуальную на заре эргодической теории  задачу поиска  
 действий  с одинаковым спектром, но метрически (от слова мера) неизоморфных. Вскоре   Анзаи \cite{Anz} вводит понятие косого произведения преобразований и, в частности, находит  решение этой задачи.  Его  реализует  косое произведение   $T$, не изоморфное  $T^{-1}$,
так как коммутативное действие и его обратное  спектрально изоморфны. 
 Познее Колмогоров  \cite{Ko} совершает энтропийный прорыв в теории  динамических систем,  получая 
континуум метрически  неизоморфных преобразований с одинаковым лебеговским
спектром.  Энтропия является симметричным числовым инвариантом, так как ее значения для  преобразований $T$ и  $T^{-1}$  одинаковы. Отметим, что главной  целью настоящей заметки является поиск несимметричного энтропийного инварианта. Для $\Z$-действий в отличие от действий групп $\Z^n$ при $n>1$,  эта задача пока не решена. Имеется  надежда  собрать асимметрию по крупицам подобно тому, как  это   было сделано в \cite{aa} для других целей. Нужно найти энтропийный инвариант такой, что для подходящего действия энтропийные характеристики   вдоль последовательности $\P$ постепенно убывали к нулю, т.е. $\P$-инвариянт  был бы нулевой, а для симметричной последовательности $\P^\ast=-\P$ он был бы положительным. 

Перемешивающие $\Z$-действия, не обладающие кратным перемешиванием, до сих пор не найдены (о проблеме Рохлина и ее взамосвязях с другими инвариантами динамических систем см.  \cite{Ry}).    Разнообразные свойства  энтропийного типа  могли бы появиться до введения Колмогоровым энтропии действия  как усиление   свойства  кратного перемешивания, если   позволить    кратности перемешивания медленно расти.  Инварианты, играющие  подобную роль, возникли  в работах  Кириллова \cite{Ki} и Кушниренко \cite{Ku}  с определенным опозданием как некоторое обобщение и реакция на ставшее  фундаментальным  понятие  энтропии на единицу времени.   Кириллов дал определение для действий групп, а Кушниренко обнаружил конечность $2^j$-энтропии автоморфизмов  $T$, входящих а орициклическтй поток.  Тем самым  в классе систем с классической нулевой энтропией им были   обнаружены  неизоморфные системы с одинаковым лебеговским спектром.
Таковыми являются степени $T\times\dots\times T$, так как значения $2^n$-энтропии у них различны.

Другое важное событие, связанное с проблемой Рохлина, энтропией и  асимметрией  действия состоялось еще позже. Идея Рохлина использовать автоморфизмы коммутативных групп для решения проблемы о  кратном  перемешивании оказалась    успешной  в случае $\Z^n$-действий  при $n>1$. Соответствующие примеры обнаружил Ледрапье. В его статье  \cite{L} указан элегантный  пример симметричного действия, но очевидным образом с этим примером  ассоциирован бесконечный  класс систем, включающий несимметричные действия. 
 Все эти действия  называем системами   Ледрапье. Приведем наиболее простой пример, который  будет использован нами в дальнейшем. 

Рассмотрим   группу $X\subset \Z_2^{\Z^2}$, состоящую из  всех последовательностей $x$, удовлетворяющих тождеству  
$$x+T^{(1,0)}x + T^{(0,1)}x =0. \eqno (1) $$
Действие $\{T^z\}$  задано   равенством $T^zx(v)=x(v-z)$.  Оно  сохраняет меру Хаара $\mu$, а
групповые автоморфизмы $T^z$  тем самым становятся автоморфизмами пространства $(X,\mu)$.
Действие  несимметрично в том смысле, что для его обратного  тождество  (1) нарушается.  Мы воспользуемся этим и установим, что для специально выбранной последовательности $\P$ конечных подмножеств группы $\Z^2$ наше действие обладает нулевой $\P$-энтропией,
а его обратное обладает    бесконечной вполне положительной $\P$-энтропией. Из этого вытекает не только неизоморфизм, но более сильное свойство дизъюнктности действия и его обратного.

Итак, кратное перемешивание и его отсутствие, энтропия и асимметрия действия в наиболее сильном ее проявлении преплетаются в классе систем Рохлина-Ледрапье.  К сказанному мы добавим типичность в смысле Альперна-Тихонова подобных эфектов.  Мы  покажем, что  
асимметрия $\P$-энтропии  для некоторых последовательностей $\P$  типична в пространстве перемешивающих  $\Z^n$-действий при $n>1$. Отсюда вытекает для  $n>0$ типичность  дизъюнктности $\Z^n$-действия и его  обратного. Здесь применяются   результаты Баштанова \cite{B} и Тихонова   \cite{T}, наша теорема об асимметрии энтропии и   $\P$-аналог теоремы Пинскера  из работы Рыжикова-Тувено, \cite{RT}, теорема 4.1.
В случае $n=1$ для доказательств типичности сильной асимметрии действия используем пример Рудольфа \cite{Ru} перемешивающего автоморфизма, дизъюнктного со своими степенями. В литературе, насколько нам известно,
это единственный имеющийся пример.  Возможно, подходящими системами  также являются  пуассоновские надстройки  из работы \cite{aa}, но это не доказано. 

\vspace{5mm}
{{\bf \Large  2. Асимметрия  $\P$-энтропии} 

\vspace{5mm}
 Пусть заданы  последовательность $\P=\{P_j\}$ конечных множеств $P_j\subset G$, где    $G$  -- некоторая бесконечная группа, и   сохраняющее меру действие $\Phi$ 
группы $G$.  
Для конечного измеримого разбиения $\xi$  пространства $X$ 
определим величины  $$h_j(\Phi,\xi)=\frac 1 {|P_j|} H\left(\bigvee_{g\in P_j}T_z\xi\right),$$
где $H(\xi)=-\sum_{i=1}^p \mu( C_i)\ln \mu( C_i).$ — энтропия разбиения $\xi=\{C_1,\dots,C_p\}$. 

\vspace{3mm}
\bf $\P$-энтропия действия. \rm    $\Phi$ относительно разбиения $\xi$ есть величина 
$h_{P}(G,\xi)={\limsup_j} \ h_k(\Phi,\xi).$
  $\P$-энтропия  действия  $\Phi$ определяется как  
$h_{P}(\Phi)=\sup_\xi h_{P}(\Phi,\xi),$
верхняя грань берется по конечным $\mu$-измеримым разбиениям  $\xi$.  

\vspace{3mm}
\bf Примеры действий Ледрапье. \rm Мы уже описали выше $\Z^2$-действие. Рассмотрим случай $n=3$. Группа $\Z_2^{\Z^3}$ содержит  компактную  подгруппу 
$X$, образованную всеми последовательностями $x: \Z^3\to\Z_2$, удовлетворяющими  тождеству
 $$x(z+(1,0,0))+x(z+(0,1,0)) +x(z+(0,0,1)) =x(z).$$
Действие группы $\Z^3$ определено равенством 
$T^wx(z)=x(z-w).$ Несложно проверить выполнение индуцированных тождеств
$$x(z+(2^n,0,0))+x(z+(0,2^n,0)) +x(z+(0,0,2^n)) =x(z).$$
Это замечательное 2-подобие играет ключевую роль для дальнейшего.
Аналогичным образом определяются  нужные нам $\Z^n$-действия  при $n>3$.
Далее мы будем рассматривать только $\Z^2$-действия Ледрапье, так как другие случаи    аналогичны, но формулы длиннее.  

\vspace{3mm}
\bf Выбор последовательности $\P$. \rm
Определим $P_j\subset \Z^2$  по индукции:
$$P_0=\{(0,0), (0,1), (1,0)\},$$
$$P_{j+1}=  2^{j}P_j\ \cup \ 2^{j}P_j+( 4^{j}, 0)\  \cup \ 2^{j}P_j+(0, 4^{j}).$$ 
Пусть $\xi=\{A,\, X\setminus A\}$, где  $A=\{x\in X\, :\, x(0,0)=0\}.$
Обозначим $P_j^\ast=-P_j$.

\vspace{3mm}
\bf Теорема 2.1. \it  $\P$-энтропия действия $\Phi$ нулевая, а $\P$-энтропия действия $\Phi^\ast$ вполне положительна.\rm

\vspace{3mm}
Теорема непосредственно вытекает из следующих лемм.

\vspace{3mm}
\bf Лемма 2.2. \it Пусть разбиение $\nu$ подчинено конечной сумме сдвигов разбиения $\xi$.  Для всех достаточно больших $j$ выполняется 
$Н(\vee_{z\in P^\ast_j} T_z \nu)=3^j H(\nu)$.\rm

\vspace{3mm}
В пространстве $L_2(X,\mu)$ характероы образуют  
ортонормированный  базис. Рассмотрим  характеры  $\chi_z$  группы $X$, определенные равенствами 
$$\chi_z(x)=1, \ x(z)=0; \ \ \ \chi_z(x)=-1, \ x(z)=1.$$
Пусть $\xi=\{A,\, X\setminus A\}$, где  $A=\{x\in X\, :\, x(0,0)=0\}.$
Группе характеров отвечает минимальная алгебра множеств, относительно которой эти характеры измеримы. Эта алгебра  порождается атомами  одинаковой меры,
 что легко доказывается по индукции. Надо лишь учесть, что разные характеры ортогональны и то, что  что все характеры вида 
$\chi_{(q,0)}$  независимы: всевозможные  их произведения являются нетривильными характерами. Алгебра множеств, ассоциированная с группой, порожденной  характерами $\chi_z$, $z\in P^\ast_j$, 
содеожит $3^j$  атомов, которые суть атомы разбиния $\xi_j=\bigvee_{z\in P_j} T^z\xi$.  Cдвиг алгебры под действием автоморфизма, дает аналогичный сдвиг характеров. Но множество $P_j$  и его сдвиги далеко расположены друг от друга при больших $j$.   Группа характеров, ее сдвиг по горизонтали и сдвиг по вертикали   независимы, поэтому   независимы соответствующие разбиения  
$$\xi_j, \ \  T^{(0,4^j)}\xi_j, \ \   T^{(0,4^j)}\xi_j.$$
Все это приводит к утверждению леммы (детали оставляем заинтересованному читателю).

\vspace{3mm}
\bf Лемма 2.3. \it  $h_j(\Phi,\xi)= (2/3)^j \ln 2.$
\rm 

\vspace{3mm}
Доказательство. 
 Заметим, что группа $G_{P_j}$, порожденная по определению всеми  характерами $\chi_z$ при   $z\in P_j$, изоморфна группе  
$G_{Q_j}$, где $Q_j= P_j\cap  \{ (q,0): q\in \Z 
\}$, так как $$\chi_{(q,4^m)}= \chi_{(q,0)}\chi_{(q,4^m)}.$$ 
Все характеры из нашей группы над линией $D=\{(w,0): w\in \Z\}$ суть произведения характеров вида $\chi_d, d\in D$.

  Число атомов разбиения $\vee_{z\in P^\ast_j} T_z \xi$ равно $2^n$, 
следовательно, действие имеет нулевую $\P$-энтропию.

\newpage
{\bf \Large 3. \ Типичноcть асимметрии }

\vspace{5mm}
Ряд свойств  перемешивающих   $\Z^n$-действий  групп $\Z^n$  на  вероятностном  пространстве $(X,\B,\mu)$ является типичным. Таковыми являются свойство  кратного перемешивания, что установил Тихонов (см. \cite{T} и ссылки), и например,   тривиальность централизатора действия, что показал  Баштанов \cite {B}.   Типичность свойства означает, что все действия с этим свойством содержат плотное множество типа $G_\delta$.  Топология в пространстве  перемешивающих $\Z^n$-действий   задается определяемой ниже метрикой $r_n$. 
  
Группа  $Aut(\mu)$  атоморфизмов пространства $(X,\mu)$ оснащена метрикой Халмоша $\rho$: для $S,T$: 
$$ \rho(S,T)=\sum_i 2^{-i}\left(\mu(S^A_i\Delta TA_i)+\mu(S^{-1}A_i\Delta T^{-1}A_i)\right),$$
где $\{A_i\}$  --   фиксированное множество, плотное в алгебре $\B$.
Складывая соответственно пастояния между образующими $\Z^n$-действий,
мы получаем метрику $\rho_n$  на пространстве $\Z^n$-действий.

Введем  на $Aut(\mu)$ метрику  $d$:
$$ d(S,T)=\sum_{i,j} 2^{-i-j}\left|\mu(SA_i\cap A_j)-\mu(TA_i\cap A_j)\right|.$$
На множестве $Mix$ всех перемешивающих  $\Z^n$-действий метрика $r_n$ задается следующим образом:
$$  r_n(\{S^z\},\{T^z\})= \rho_n(\{S^z\},\{T^z\}) + \sup_{z\in \Z^n}d(S^z,T^z).$$
Напомним, что  действие  $\Phi=\{T^z\,:\, z\in\Z^n\}$ называется перемешивающим, если для всяких $A,B\in\B$ выполнено
$$ \mu(T^zA\cap B)\to \mu(A)\mu(B), \ z\to\infty.$$

Для  перемешивающих $\Z^n$-действий известна (и мы ее заново докажем в частном случае)    типичность  бесконечной энтропии Кириллова-Кушниренко   для последовательностей расширяющися множеств $P_j$ 

Рассмотрим   пространство $\M$  мер на $X\times X$,
удовлетворяющих следующим проекционным условиям:  для всякого $A\in \B$ выполнено
$$\nu(A\times X)=\nu(X\times A)=\mu(A).$$
Множество   $\M$  оснащается метрикой  $dist$:
$$ dist(\nu,\eta)=\sum_{i,j} 2^{-i-j}\left|\nu(A_i\times A_j)-\eta(A_i\times A_j)
\right|,$$  где $\{A_i\}$  -- некоторое  фиксированное семейство,  плотное в алгебре всех $\mu$-измеримых множеств.  Пространство $(\M,dist)$ является компактом. 

\it Дизъюнктность \rm автоморфизмов $S,T$  означает, что для всех мер, кроме $\mu^2=\mu\times \mu$,  выполнено   $dist(\nu, (S\times T)\nu)>0$.

\vspace{3mm}
\bf Теорема 3.1.  \it Типичное перемешивающее    $\Z^n$-действие  и его обратное  дизъюнктны.\rm

\vspace{3mm}
Доказательство.  Теорема при $n>1$ является следствием теоремы 2.1 и теоремы 4.1 из \cite{RT} и следующих фактов. Легко убедиться в том, что  действия со свойствами, которые  утверждает теорема  2.1, образуют $G_\delta$-множество.   Класс сопряженности перемешивающего действия плотен в пространстве перемешивающих действий, как показал Баштанов \cite{B}.  Можно также воспользоваться результатами  Тихонова \cite{T}, среди которых  плотность класса сопряженности декартового квдрата  перемешивающего $\Z^n$-действия (и $\R^n$-действия).
Квадрат  наследует нулевую энтропию  действия при однвременной вполне положительной  энтропии обратного действия. Это обеспечивает плотность уазанного  $G_\delta$-множества. Отметим, что этот подход  эффективен  в доказательстве аналога  теоремы 3.1 для $\R^n$-действий (для которых нет теоремы Баштанова).

 Осталось рассмотреть случай $n=1$. Мы докажем следующее утверждение и тем самым теорему  3.1. 

\vspace{3mm}
\bf Теорема 3.2.  \it Для типичного перемешивающего автоморфизма $T$ вероятностного
пространства его степени  $T^p$ и $T^q$ при  $p\neq q$,  дизъюнктны.\rm

\vspace{3mm}
Доказательство. 
Пример преобразования  с указанным свойством  дан в работе Рудольфа \cite{Ru}.
Класс сопряженности всякого перемешивающего автоморфизма плотен в $(Mix, r_1)$, как показал Баштанов \cite{B}.  Убедимся в том, что все автоморфизмы с указанным свойством  образуют $G_\delta$-множество.   Рассмотрим только случай $p=1, q=-1$, так как доказательство остальных  случаев совершенно аналогично.    Положим
$$K_m=\{\nu\in \M\ :\ dist(\nu,\mu^2)\ \geq \ 1/m\}.$$
Множества $K_m$ замкнуты  в $(\M,dist)$, следовательно, они компактны. Для фиксированного $T\in Aut$,
дизъюнктного с $T^{-1}$,   выполнено
$$ \min \{  dist(\nu\, ,\, (T\times T^{-1})\nu) \ :\ \nu\in K_m\}=d_m>0,$$
так как непрерывная положительная функция на компакте $K_m$  отделена от 0.
Определим  множество 
$$U_{m}= \bigcup_n  \{ S\,:\,  dist(\nu, \, (S\times S^{-1})\nu)> \ 1/n , \ \nu\in K_m\}.$$
Оно открыто как объединение открытых множеств.
Остается заметить, что  $\bigcap_m U_m$  есть семейство   всех перемешивающих автоморфизмов,  дизъюнктных со своими обратными, является  плотным $G_\delta$-множеством.  Теорема доказана.

\vspace{2mm}
\bf Теорема 3.3.  \it Пусть  $a$ -- эндоморфизм группы $\Z^n$. Если найдется перемешивающее действие $\Phi\in Mix_n$, дизъюнктное с   $\tilde\Phi^a$, то   множество всех перемешивающих действий $\Phi$, дизъюнктных с   $\Phi^a$, типично.\rm

\vspace{2mm}
Это условное утверждение  доказывается тем же методом, что и теорема 3.1. Плотность
класса сопряженности перемешивающего  $\Z^n$-действия установлена в \cite{B}.
Но нужные  примеры действий в  литературе отсутствуют.  Отметим, что   их можно  построить,  но это  требует значительных усилий.  Возможно действия Ледрапье  обеспечивают 
теорему нужными примерами.

\section{ Замечания}
 Отметим маленький квази-парадокс. Мы получили асимметрию для типичных действий из отсутствия   
кратного перемешивания. Но типичные  перемешивающие $\Z^n$-действия, как  показал Тихонов \cite{T12},  обладают перемешиванием всех кратностей.  

Объясним  название нашей статьи. Из теоремы  3.1.1 \cite{R21} (в  формулировке этой теоремы случайно было опущено слово "аменабельная")  
  вытекает, что  для типичного (относительно  метрики Халмоша)    $\Z^n$-действия  
$\{P_j\}$-энтропия и 
$\{P_j^\ast\}$-энтропия обе бесконечны.  Для перемешивающих $\Z^n$-действий доказываем  и применям теорему Баштанова о всюду плотностие  класса сопряженности.   В чем же проявляется асимметрия?  Ответ дает теорема, которую мы концептуально  доказали, поэтому ограничимся только формулировкой.

\vspace{3mm}
\bf Теорема 4.1. \it Для нашей  последовательности $\P=\{P_j\}$ типично множество  перемешивающих действий  $\Phi$, для которых найдется собственная  последовательность $j(k)\to\infty$ такая, что   $\{P_{j(k)}\}$-энтропия действия $\Phi$ нулевая, а  $\{P^\ast_{j(k)}\}$-энтропия  вполне положительна. В этом случае у обратного действия $\Phi^\ast$ 
все наооборот:  $\{P^\ast_{j(k)}\}$-энтропия нулевая, а  $\{P_{j(k)}\}$-энтропия вполне положительна.  При этом $P_j$ и $P_j^\ast$ энтропии этих действий бесконечны.\rm
\vspace{3mm}

\vspace{3mm}
Благодаря $\P$-энтропии можно усилить эффект энтропии Колмогорова, различающей континуум   действий с лебеговским спектром. Можно указать континуум попарно дизъюнктных действий с лебеговским спектром \cite{faa}. Но эти действия симметричны, так как являются факторами всегда  симметричных гауссовских систем.


\begin{thebibliography}{9}
\bibitem{Ki}A.A. Kirillov, Dynamical systems, factors, and representations of groups, Uspekhi Mat. Nauk, 22:5(137) (1967), 67–80
\bibitem{Ku} A.G. Kushnirenko, On metric invariants of entropy type, Uspekhi Mat. Nauk, 22:5 (137) (1967), 57–65
\bibitem{RT} V.V. Ryzhikov and J.-P. Thouvenot, Quasi-similarity, entropy, and disjointness of ergodic actions, Funkts. Analysis and Its Appl., 58:1 (2024), 117–124

\bibitem{Ro} V.A. Rokhlin, On Endomorphisms of Compact Commutative Groups, Izv. Akad. Nauk SSSR, Ser. Mat., 13:4 (1949), 329–340

\bibitem{Anz} H. Anzai, Ergodic Skew Product Transformations on the Torus. Osaka Math. J., 1951, 3:1, 83–99

\bibitem{Ry} 
V.V. Ryzhikov, Multiple mixing, local rank, and joinings of probability space automorphisms, Theory Probab. Appl., 70:4 (2026), 614–622

\bibitem{Ko} A.N. Kolmogorov, A new metric invariant of transitive dynamical systems and automorphisms of Lebesgue spaces, Dokl. Akad. Sci. USSR, 119:5 (1958), 861–864

\bibitem{L} F. Ledrappier, Un champ markovien peut \^etre d'entropie nulle et 
m\'elangeant, C. R. Acad. Sci. Paris S\'er. A, 287:7 (1978), 561–563

\bibitem{aa} V.V. Ryzhikov, Asymmetric mixing Poisson suspensions, Algebra and Analysis, 37:4 (2025), 141–148

\bibitem{B} A.I. Bashtanov, Conjugacy Classes are Dense in the Space of Mixing $Z^d$-Actions, Math. Notes, 99:1 (2016), 9–23

\bibitem{T} S.V. Tikhonov, On mixing actions of locally compact groups, Mat. Notes, 117:4 (2025), 561-574


\bibitem{Ru} Rudolph D., An example of measure-preserving map with minimal self-joinings, and applications, J. d'Analyse Math., 35 (1979), 97-12

\bibitem{T12}S.V. Tikhonov, Genericity of a multiple mixing, Russian Math. Surveys, 67:4 (2012), 779-780

\bibitem{R21}V.V. Ryzhikov, Compact families and typical entropy invariants of measure-preserving actions, Trans. Moscow Math. Soc., 82 (2021), 117-123
\bibitem{faa}V.V. Ryzhikov, Continuum disjoint automorphisms with Lebesgue spectrum,
Funct. Anal. Appl., 60 (2026)


\end{thebibliography}

\newpage
\end{document}